\documentclass[12pt]{article}
\usepackage[T2A]{fontenc}
\usepackage[utf8]{inputenc}
\usepackage{amssymb,amsmath,amsfonts,amsthm,amscd,latexsym,indentfirst}
\usepackage{epsfig,geometry}
\usepackage{wasysym,verbatim,lineno}
\usepackage[all]{xy}
\def\RB{\textrm{RB}}

\def\Var{\mathrm{Var}}
\def\Span{\mathrm{Span}}

\def\pre{\mathrm{pre}}
\def\post{\mathrm{post}}
\def\As{\mathrm{As}}
\def\Lie{\mathrm{Lie}}
\def\Com{\mathrm{Com}}

\begin{document}

\begin{flushright}
\end{flushright}

\begin{center}
{\Large
Poincar\'{e}---Birkhoff---Witt theorem \\ for pre-Lie and postLie algebras}

V. Gubarev
\end{center}

\begin{abstract}
We construct the universal enveloping preassociative and postassociative algebra
for a pre-Lie and a postLie algebra respectively.
We show that the pairs $(\pre\Lie,\pre\As)$ and $(\post\Lie,\post\As)$ 
are Poincar\'{e}---Birkhoff---Witt-pairs, for the first one it's a reproof of the result of V. Dotsenko and P. Tamaroff.

\end{abstract}

\section{Introduction}

In 1960s, pre-Lie algebras appeared independently in affine geometry
(E. Vinberg~\cite{Vinberg63}; J.-L. Koszul~\cite{Koszul61}),
and ring theory (M. Gerstenhaber~\cite{Gerst63}).
Arising from diverse areas, pre-Lie algebras are known under
different names like Vinberg algebras, Koszul algebras,
left- or right-symmetric algebras (LSAs or RSAs), Gerstenhaber algebras.
Pre-Lie algebras satisfy an identity
$(x_1 x_2)x_3 - x_1(x_2 x_3) = (x_2 x_1) x_3 - x_2 (x_1 x_3)$.
See \cite{Burde06,Manchon} for surveys on pre-Lie algebras.

In 2001, J.-L. Loday \cite{Dialg99} defined the dendriform (di)algebra
(preassociative algebra) as a vector space endowed 
with two bilinear operations $\succ,\prec$ satisfying 
$$
\begin{gathered}
(x_1\succ x_2+x_1\prec x_2)\succ x_3 = x_1\succ (x_2 \succ x_3), \
(x_1\succ x_2)\prec x_3=x_1\succ(x_2\prec x_3), \\
x_1\prec(x_2\succ x_3+x_2\prec x_3)=(x_1\prec x_2)\prec x_3.
\end{gathered}
$$

In 1995, J.-L. Loday also defined~\cite{Loday95} Zinbiel algebra
(precommutative algebra), on which the identity
$(x_1x_2 + x_2x_1)x_3 = x_1(x_2 x_3)$ holds. 
Every preassociative algebra with the identity $x\succ y = y\prec x$ 
is a precommutative algebra ($x_1x_2 = x_1\succ x_2$)
and under the product $x\cdot y = x\succ y - y\prec x$ is a pre-Lie algebra.

In 2004, dendriform trialgebra (postassociative algebra) was introduced~\cite{Trialg01}, 
i.e., an algebra with bilinear operations $\prec,\succ,\cdot$ satisfying seven certain axioms.
A space~$A$ with two bilinear operations $[,]$ and $\cdot$ is called a post-Lie algebra 
(B.~Vallette, 2007 \cite{Vallette2007}) if $[,]$ is a Lie bracket and the next identities hold
$$
(x \cdot y) \cdot z - x \cdot (y \cdot z)
- (y \cdot x) \cdot z + y \cdot (x \cdot z) = [y,x]\cdot z, \quad
x\cdot [y,z] = [x \cdot y,z] + [y,x\cdot z].
$$
In last dozen years, an amount of articles devoted to post-Lie algebras
in different areas is arisen \cite{Burde12,ELMM,sl2}. 

Let us explain the choice of terminology. 
Given a binary quadratic operad~$\mathcal P$, 
the defining identities for pre- and post-$\mathcal P$-algebras were found in \cite{BBGN2012}.
One can define the operad of pre- and post-$\mathcal P$-algebras
as $\mathcal P \bullet \pre\Lie$ and
$\mathcal P \bullet \post \Lie$ respectively.
Here $\pre\Lie$ and $\post\Lie$ denote 
the operads (varieties) of pre-Lie algebras and postLie algebras respectively,
$\bullet$ denotes the black Manin product of operads \cite{GinzKapr}.
By pre- or postalgebra we will mean pre- or post-$\mathcal P$-algebra 
for some operad~$\mathcal P$.

Before stating the main problem of the work we introduce a very useful tool
to deal with pre- and postalgebras, so called Rot---Baxter operators, and 
the notion of a~Poincar\'{e}---Birkhoff---Witt pair.

A linear operator $R$ defined on an algebra $A$ over a field $\Bbbk$
is called a Rota---Baxter operator (RB-operator, for short) of a weight $\lambda\in\Bbbk$
if it satisfies the relation
\begin{equation}\label{RB}
R(x)R(y) = R( R(x)y + xR(y) + \lambda xy), \quad x,y\in A.
\end{equation}
In this case, an algebra $A$ is called Rota---Baxter algebra (RB-algebra).

G.~Baxter defined the notion of what is now called Rota---Baxter operator on 
a (commutative) algebra in 1960 \cite{Baxter60},
solving an analytic problem. The relation~\eqref{RB} with $\lambda = 0$
appeared as a generalization of integration by parts formula.
G.-C.~Rota~\cite{Rota68}, P.~Cartier \cite{Cartier72} and others studied 
different combinatorial properties of RB-opera\-tors and RB-algebras.
In 1980s, the deep connection between Lie RB-algebras and Yang---Baxter equation
was found \cite{BelaDrin82,Semenov83}. More about Rota---Baxter algebras 
see in the monograph of L.~Guo~\cite{Guo2011}.

In 2000, M. Aguiar \cite{Aguiar00} stated that an associative 
algebra with a given Rota---Baxter operator~$R$
of weight zero under the operations $a \succ b = R(a)b$, $a\prec b = aR(b)$ 
is a~preassociative algebra. In 2002, K. Ebrahimi-Fard~\cite{Fard02}
showed that an associative RB-algebra of nonzero weight $\lambda$ 
under the same two products $\succ$, $\prec$ and the third operation $a\cdot b = \lambda ab$ 
is a~postassociative algebra.
The analogue of the Aguiar construction for the pair of pre-Lie algebras
and Lie RB-algebras of weight zero was stated in 2000 by M.~Aguiar \cite{Aguiar00} 
and by I.Z. Golubchik, V.V.~Sokolov~\cite{GolubchikSokolov}.
In 2010 \cite{BaiGuoNi10}, this construction for the pair of post-Lie algebras 
and Lie RB-algebras of nonzero weight was extended.

In 2013 \cite{BBGN2012}, the construction of M. Aguiar and K. Ebrahimi-Fard
was generalized for the case of arbitrary variety.

In 2008, the notion of universal enveloping RB-algebras
of pre- and postassociative algebras was introduced~\cite{FardGuo07}.
In~\cite{FardGuo07}, it was also proved that the universal enveloping of 
a free pre- or postassociative algebra is free.

In 2010, with the help of Gr\"{o}bner---Shirshov bases \cite{BokutChen},
Yu. Chen and Q.~Mo proved that every preassociative algebra
over a field of characteristic zero injectively embeds into
its universal enveloping RB-algebra \cite{Chen11}.

In 2013 \cite{GubKol2013}, given a variety $\Var$, it was proved that every pre-$\Var$-algebra
(post-$\Var$-algebra) injectively embeds into its universal enveloping $\Var$-RB-alge\-bra of weight 
$\lambda = 0$ ($\lambda\neq0$). 
Further, author constructed universal enveloping RB-algebra for a given pre- or postalgebra 
in commutative \cite{Gub2017Com}, associative \cite{Gub2017As}, and Lie \cite{Gub2017Lie} cases.
In the associative case it gave an answer to the question of L.~Guo \cite[p. 148]{Guo2011}.

The classical Poincar\'{e}---Birkhoff---Witt theorem states that
given a Lie algebra $g$ with a linear basis $\{x_i \mid i \in I\}$, where $I$ is a well-ordered set, 
the monomials $x_{i_1}\ldots x_{i_n}$ with $i_1\leq \ldots \leq i_n$
form a linear basis for the universal enveloping associative algebra $U(g)$.
As a consequence, we get that the linear basis of the algebra $U(g)$
does not depend on the product in the Lie algebra $g$.
Such relationships between two varieties $\mathcal{V},\mathcal{W}$ 
in the case when there exists a functor $\phi\colon \mathcal{V}\to \mathcal{W}$ 
associating to every algebra $A\in \mathcal{V}$ an algebra $\phi(A)\in\mathcal{W}$
by changing multiplication in $A$ was generalized by I. Shestakov and A.A.~Mikhalev 
in the term Poincar\'{e}---Birkhoff---Witt (PBW-) pair, see details in~\cite{PBW}.

Now let us formulate the main problem to which the work is devoted.
In advance, $\pre\As$ and $\post\As$
denote the varieties of pre- and postassociative algebras respectively.

{\bf Problem 1}.
a) Prove that every pre-Lie (postLie) algebra injectively embeds into
its universal enveloping preassociative (postassociative) algebra.

b) Clarify if the pairs $(\pre\As,\pre\Lie)$ and $(\post\As,\post\Lie)$ are PBW-ones.

c) Construct the universal enveloping preassociative 
(postassociative) algebra for a~given pre-Lie (postLie) algebra. 

For pre-Lie algebras, Problem~1b and special version of Problem~1c were stated 
by P.~Kolesnikov in~\cite{Kol2017} in the context of Gr\"{o}bner---Shirshov bases 
for preassociative algebras. 
J.-L. Loday asked V. Dotsenko about the solution of Problem~1b around 2009~\cite{Dots}.
The discussion of Problem~1 in the case of restricted 
pre-Lie algebras can be found in \cite{Dokas13}. 
The analogues of Problem~1 for Koszul-dual objects, di- and trialgebras, 
were solved in \cite{LP93,GubKol2014}.

Recently, in~\cite{Gub2018-pre} and \cite{Gub2018-post}, 
the author solved Problem 1a in pre- and postalgebra cases
with the help of embedding of pre-Lie (postLie) algebras into Lie RB-algebras~\cite{GubKol2013}
and the Gr\"{o}bner---Shirshov bases technique developed for associative RB-algebras~\cite{Guo2013}.
Actually, the solution of Problem~1a for pre-Lie algebras can be derived from the results concerned
Hopf preassociative algebras and so called brace algebras stated in 2002
independently by F.~Chapaton~\cite{Chapaton} and M. Ronco~\cite{Ronco}.
In 2018, V.~Dotsenko and P.~Tamaroff by means of the general approach arising from the category theory 
solved Problem~1b for pre-Lie algebras stating that the pair of varieties 
$(\mathrm{pre}\Lie,\mathrm{pre}\As)$ is a PBW-pair~\cite{Dots}.

The current work is devoted to the complete solution of Problem~1 in both pre- and postalgebra cases.
Let us briefly describe the idea of the solution. For this, we need one more embedding problem.

Let $A$ be an associative algebra with an RB-operator $R$.
Then the algebra $A^{(-)}$ is a Lie RB-algebra
under the product $[x,y] = xy-yx$ and the same action of $R$. 
Thus, we can state the analogue of Problem~1 for 
the varieties of Lie and associative RB-algebras. 

{\bf Problem 2}.
a) Prove that every Lie RB-algebra injectively embeds into
its universal enveloping associative RB-algebra.

b) Construct the universal enveloping associative RB-algebra 
for a given Lie RB-algebra. 

We are not asking whether the pair $(\RB\As,\RB\Lie)$ 
of the varieties of associative and Lie RB-algebras
forms a PBW-pair, since it is easy to disprove it just in 2-dimensional case.
To the moment, we are far from the solution Problem~2, and the current work 
as well as~\cite{Gub2018-post,Gub2018-pre} can be also considered like a step in such direction. 

The sketch of the solution of Problem~1c is following. 
At first, we embed a pre- or post-Lie algebra $C$ into its universal enveloping Lie RB-algebra $L$~\cite{Gub2017Lie}. 
At second, with the help of Gr\"{o}bner---Shirshov bases 
we embed the Lie RB-algebra $L$ into its universal enveloping associative RB-algebra $A$
with an RB-operator $R$. At third, we show that a~subalgebra $U(C)$ generated by~$C$
in the induced pre- or postassociative algebra on the space $A$ is the universal 
enveloping pre- or postassociative algebra for $C$.

Actually, the same solution algorithm was used earlier by author in~\cite{Gub2018-post} 
and~\cite{Gub2018-pre} but with the following change: in the first step 
it was considered injective enveloping Lie RB-algebra from~\cite{GubKol2013} 
instead of $L$, the universal enveloping one constructed in~\cite{Gub2017Lie}. 
Surprisingly, a posteriori we may say that at least in the pre-Lie algebra case 
the injective enveloping preassociative algebra of a given pre-Lie algebra
obtained earlier in~\cite{Gub2018-pre} is an universal one.

At the end of Introduction, let us collect all stated in the work connections between universal enveloping 
algebras of different kind in the following commutative diagram,

\medskip
\begin{xy}
(225,20)*+{\pre/\post\Lie}="a"; (255,20)*+{\RB\Lie}="b";%
(225,0)*+{\pre/\post\As}="c"; (255,0)*+{\RB\As}="d";%
{\ar "a";"b"};
{\ar "a";"c"};{\ar "b";"d"};%
{\ar "c";"d"};%
{\ar@{->} "a";"d"};
\end{xy}

\medskip 

\noindent where every arrow maps the algebra from corresponding variety into its universal enveloping one.
An associative RB-algebra $A$ with an RB-operator $P$ of weight $\lambda = 0$ ($\lambda \neq 0$) 
is an enveloping for a~pre-Lie (postLie) algebra $C$ in the sense that the following equalities 
\begin{equation}\label{pre/postLie2RBAs}
a\succ b = P(a)b - bP(a), \quad
a\prec b = aP(b) - P(b)a \quad 
(a\cdot b = \lambda ab - \lambda ba)
\end{equation}
hold for any $a,b\in C$. 

\section{Preliminaries}

\subsection{Some required formulas}

Given a Lie algebra $L$, denote the product $[[\ldots[[y,x],x]\ldots],x]\in L^{p+1}$ by $[y,x^{(p)}]$.

{\bf Lemma}~\cite{Gub2018-pre}.
Given a Lie algebra~$L$, the equality
\begin{equation}\label{Lemma}
(l+1)yx^l
 = \sum\limits_{i=2}^{l+1}(-1)^i \binom{l+1}{i}[y,x^{(i-1)}]x^{l+1-i}
 + \big(yx^l+xyx^{l-1} + \ldots + x^l y)
\end{equation}
holds in the universal enveloping algebra $U(L)$ for any $x,y\in L$ and $l\geq0$.

It follows immediately from~\eqref{RB} that
\begin{multline}\label{LongRB}
R(b_1)\ldots R(b_t) 
 - R\bigg( \sum\limits_{i=1}^t R(b_1)\ldots R(b_{i-1})\hat{R}(b_i)R(b_{i+1})\ldots R(b_t) \\
 + \lambda \sum\limits_{1\leq i_1<i_2\leq t} R(b_1)\ldots \hat{R}(b_i)\ldots \hat{R}(b_j)\ldots R(b_t) + \ldots \\
 + \lambda^{k-1}\sum\limits_{1\leq i_1<i_2<\ldots<i_k\leq t} R(b_1)\ldots \hat{R}(b_{i_1})\ldots \hat{R}(b_{i_k})\ldots R(b_t) 
 + \ldots +\lambda^{t-1}b_1\ldots b_t \bigg) = 0,
\end{multline}
where the sign $\hat{}$ shows the omitting action of $R$. 
From~\eqref{Lemma} and~\eqref{LongRB}, we derive the formula
\begin{multline}\label{LongRB-conseq2}
R(R(b_1)R(b_2)\ldots R(b_{l-1})b_l R(b_l)^k R(b_{l+k+1})) \\ \allowdisplaybreaks
 = \frac{1}{k+1} R(b_1)R(b_2)\ldots R(b_{l-1})R(b_l)^{k+1} R(b_{l+k+1}) \\ 
 + \frac{1}{k+1} R\bigg(\sum\limits_{i=2}^{k+1}(-1)^i\binom{k+1}{i} R(b_1)R(b_2)\ldots R(b_{l-1})[b_l,R(b_l)^{(i-1)}]R(b_l)^{k+1-i} R(b_{l+k+1})  \\ 
 - \sum\limits_{1\leq i\leq l-1,\,i=l+k+1} R(b_1)\ldots R(b_{i-1})\hat{R}(b_i)R(b_{i+1})\ldots R(b_t) \\ 
 - \lambda \sum\limits_{1\leq i_1<i_2\leq l+k+1} R(b_1)\ldots \hat{R}(b_i)\ldots \hat{R}(b_j)\ldots R(b_{l+k+1}) - \ldots \\ 
 - \lambda^{s-1}\sum\limits_{1\leq i_1<\ldots<i_s\leq l+k+1} R(b_1)\ldots \hat{R}(b_{i_1})\ldots \hat{R}(b_{i_s})\ldots R(b_{l+k+1}) \\ 
 - \ldots - \lambda^{l+k}b_1\ldots b_{l+k+1} \bigg),
\end{multline}
where $b_l = b_{l+1} = \ldots = b_{l+k}$
and $R(b_l)^s$ means $(R(b_l))^s$. 

In advance we will use the formula~\eqref{LongRB-conseq2}
with maybe absent $R(b_1)$ and $R(b_{l+k+1})$, it means that 
we omit all $b_1$ and $R(b_1)$ in the summands and all indexes $i_j$ start with two
when $R(b_1)$ is absent. We do the same if $R(b_{l+k+1})$ is absent.

\subsection{Embedding of pre- and postalgebras into RB-algebras}

{\bf Theorem 1}~\cite{Aguiar00,BBGN2012,BaiGuoNi10,Fard02,GolubchikSokolov,Loday2007}.
Let $A$ be an RB-algebra of a variety $\Var$ and weight $\lambda=0$ ($\lambda\neq0$).
With respect to the operations
\begin{equation}\label{LodayToRB}
x\succ y = R(x)y,\quad x\prec y = xR(y)\quad (x\cdot y = \lambda xy)
\end{equation}
$A$ is a pre-$\Var$-algebra (post-$\Var$-algebra).

Denote the pre- and post-$\Var$-algebra obtained in Theorem~1 as $A^{(R)}_{\lambda}$. 

Given a pre-$\Var$-algebra $\langle C,\succ,\prec\rangle$, universal enveloping 
RB-$\Var$-algebra $U$ of $C$ is the universal algebra in the class of all RB-$\Var$-algebras
of weight zero such that there exists homomorphism from $C$ to $U^{(R)}_0$.
Analogously universal enveloping RB-$\Var$-algebra of a post-$\Var$-algebra is defined. 

{\bf Theorem 2}~\cite{GubKol2013}.
Every pre-$\Var$-algebra (post-$\Var$-algebra) could be embedded 
into its universal enveloping RB-algebra of the variety $\Var$ 
and weight $\lambda = 0$ ($\lambda\neq0$).

Let us briefly describe the idea of the proof of Theorem~2.
Given a pre-$\Var$-algebra (post-$\Var$-algebra) $C$, we define the product on the space
$\hat{C} = C\oplus C'$, where $C'$ is a copy of $C$, in such a way that 
$\hat{C}$ is an algebra of the variety $\Var$. Then we define on $\hat{C}$ 
the linear operator $P$ which occurs to be a Rota---Baxter operator 
of weight $\lambda = 0$ ($\lambda\neq0$). Finally, Theorem~2 was stated
by embedding $C$ into $\hat{C}^{(P)}_\lambda$ by the map $c\to c'$.  

\subsection{Gr\"{o}bner---Shirshov bases for associative RB-algebras}

Let $R\As\langle X\rangle$ denote the free associative algebra 
generated by a set $X$ with a linear map $R$ in the signature.
One can construct a linear basis of $R\As\langle X\rangle$ (see, e.g., \cite{FardGuo07})
by induction. At first, all elements from $S(X)$, the free semigroup generated by~$X$, 
lie in the basis. At second, if we have basic elements $a_1,a_2,\ldots,a_k$, $k\geq1$, 
then the word $w_1R(a_1)w_2\ldots w_kR(a_k)w_{k+1}$ lies in the basis 
of~$R\As\langle X\rangle$.
Here $w_2,\ldots,w_k\in S(X)$ and $w_1,w_{k+1}\in S(X)\cup \emptyset$,
where $\emptyset$ denotes the empty word.
Let us denote the basis obtained as $RS(X)$.
Given a word~$u$ from $RS(X)$, the number of appearances of the symbol~$R$
in~$u$ is denoted by $\deg_R(u)$, the $R$-degree of~$u$. 
We call an element from $RS(X)$ of the form $R(w)$ as $R$-letter.
By~$X_\infty$ we denote the union of $X$ and the set of all $R$-letters.
Given $u\in RS(X)$, define $\deg u$ (degree of $u$) as the length of $u$
in the alphabet~$X_\infty$. In~\cite{FardGuo07}, $\deg u$ was called the breadth of $u$.

Suppose that $X$ is a well-ordered set with respect to $<$.
Let us introduce by induction the deg-lex order on $S(X)$. 
At first, we compare two words $u$ and $v$ by the length:
$u < v$ if $|u|<|v|$. At second, when $|u| = |v|$, 
$u = x_i u'$, $v = x_j v'$, $x_i, x_j\in X$,
we have $u < v$ if either $x_i < x_j$ or $x_i = x_j$, $u'< v'$.
We compare two words $u$ and $v$ from $RS(X)$ 
by $R$-degree: $u<v$ if $\deg_R(u)<\deg_R(v)$. 
If $\deg_R(u) = \deg_R(v)$, we compare $u$ and $v$ in deg-lex order as words 
in the alphabet $X_\infty$. Here we define each $x$ from $X$ to 
be less than all $R$-letters and $R(a)<R(b)$ if and only if $a<b$.

Let $*$ be a symbol not containing in $X$.
By a $*$-bracketed word on $X$, we mean a basic word 
from $R\As\langle X\cup\{*\}\rangle$ with exactly one occurrence of $*$. 
The set of all $*$-bracketed words on $X$ is denoted by $RS^*(X)$.
For $q\in RS^*(X)$ and $u\in R\As\langle X\rangle$, we define
$q|_u$ as the bracketed word obtained by replacing the letter $*$ in $q$ by $u$.

The order defined above is {\it monomial}, i.e., from $u < v$  
it follows that $q|_u < q|_v$ for all $u,v\in RS(X)$ and $q\in RS^*(X)$.

Given $f\in R\As\langle X\rangle$, by $\bar{f}$ we mean the leading word in $f$.
We call $f$ monic if the coefficient of $\bar{f}$ in $f$ is~1. 

{\bf Definition 1}~\cite{Guo2013}.
Let $f,g\in R\As\langle X\rangle$.
If there exist $\mu,\nu,w\in RS(X)$ such that 
$w = \bar{f} \mu = \nu \bar{g}$ with $\deg w<\deg(\bar{f})+\deg(\bar{g})$, 
then we define $(f,g)_w$ as $f\mu - \nu g$ and call it 
the {\it composition of intersection} of $f$ and $g$ with respect to $(\mu,\nu)$.
If there exist $q\in RS^*(X)$ and $w\in RS(X)$ 
such that $w = \bar{f} = q|_{\bar{g}}$, then we define
$(f,g)^q_w$ as $f-q|_g$ and call it the {\it composition of inclusion} 
of $f$ and $g$ with respect to $q$.

{\bf Definition 2}~\cite{Guo2013}.
Let $S$ be a subset of monic elements from $R\As\langle X\rangle$ 
and $w\in RS(X)$. 

(1) For $u,v\in R\As\langle X\rangle$, we call $u$ and $v$ 
congruent modulo $(S, w)$ and denote this by $u \equiv v \mod (S, w)$ 
if $u - v = \sum c_i q_i|_{s_i}$ 
with $c_i \in \Bbbk$, $q_i\in RS^*(X)$, $s_i\in S$ and $q_i|_{\overline{s_i}} < w$.

(2) For $f,g\in R\As\langle X\rangle$ and suitable 
$w,\mu,\nu$ or $q$ that give a composition of intersection
$(f,g)_w$ or a composition of inclusion $(f,g)^q_w$, 
the composition is called trivial modulo $(S, w)$ if
$(f,g)_w$ or $(f,g)^q_w \equiv 0 \mod (S, w)$.

(3) The set $S\subset R\As\langle X\rangle$ is 
called a {\it Gr\"{o}bner---Shirshov basis} if, for all $f,g\in S$, 
all compositions of intersection $(f,g)_w$ 
and all compositions of inclusion $(f,g)^q_w$ are trivial modulo $(S, w)$.

{\bf Theorem 3}~\cite{Guo2013}. 
Let $S$ be a set of monic elements in $R\As\langle X\rangle$, 
let $<$ be a monomial ordering on $RS(X)$ and let $Id(S)$ 
be the $R$-ideal of $R\As\langle X\rangle$ generated by $S$. 
If $S$ is a Gr\"{o}bner---Shirshov basis in $R\As\langle X\rangle$,
then $R\As\langle X\rangle = \Bbbk Irr(S)\oplus Id(S)$ 
where $Irr(S) = RS(X)\setminus \{q|_{\bar{s}}\mid q\in RS^*(X),s\in S\}$
and $Irr(S)$ is a linear basis of $R\As\langle X\rangle/Id(S)$.

\section{PBW-theorem for pre-Lie and postLie algebras}

Let $A$ be an associative algebra with an RB-operator $R$.
Then the algebra $A^{(-)}$ is a Lie RB-algebra
under the product $[x,y] = xy-yx$ and the same action of $R$. 

Let $L$ be a Lie RB-algebra with an RB-operator $P$ of weight~$\lambda$.
Suppose that there exists a subset $X_0 = \{x_\alpha\mid \alpha\in \Omega\}$ in~$L$ such that
$X = \{x_{\alpha,k}:= P^k(x_\alpha)\mid k\in\mathbb{N},\,\alpha\in\Omega\}$ is a linear basis of $L$.
Our goal is to construct the universal enveloping associative RB-algebra of $L$
(via Gr\"{o}bner---Shirshov bases). This will lead us to the proof of 
the Poincar\'{e}---Birkhoff---Witt (PBW) theorem for the pairs 
(pre-Lie, preAs) and (postLie, postAs).

We may assume that the set $\Omega$ is well-ordered, so we define an order $<$ on the set $X$:
$x_{\alpha,k}<x_{\beta,l}$ if $\alpha<\beta$ or $\alpha = \beta$ and $k<l$.

Consider the set $S$ of the following elements in $R\As\langle X\rangle$:
\begin{gather}
xy - yx - [x,y],\ x>y,\ x,y\in X, \label{UnivRel} \\
R(a)R(b) - R(R(a)b + aR(b) + \lambda ab), \label{almostRB} \\
R(R(z_1)\vec{x}_{q_1}R(z_2)\ldots R(z_s)\vec{x}_{q_s}x_{\beta,r}x_{\beta,r+1}^k R(z_{s+1})) 
 - \Delta, \label{LongRel1} 
\end{gather}
where $\Delta$ is defined as the RHS of~\eqref{LongRB-conseq2} for $l=l_1+\ldots+l_s+s+1$ and
\begin{gather*}
b_1 = z_1,\
b_2 = \widetilde{x_{q_1}^1},\
b_3 = \widetilde{x_{q_1}^2},\ \ldots,\
b_{l_1+1} = \widetilde{x_{q_1}^{l_1}},\
b_{l_1+2} = z_2,\
b_{l_1+3} = \widetilde{x_{q_2}^1},\
\ldots, \\ \allowdisplaybreaks
b_{l_1+\ldots +l_{s-1}+s} = z_s,\
b_{l_1+\ldots +l_{s-1}+s+1} = \widetilde{x_{q_l}^1},\
\ldots,\ 
b_{l_1+\ldots +l_s+s} = \widetilde{x_{q_s}^{l_s}},\\ 
b_{l_1+\ldots +l_s+s+1} = \ldots 
 = b_{l_1+\ldots +l_s+s+k+1} = x_{\beta,r}, \
b_{l_1+\ldots +l_s+s+k+2} = z_{s+1}.
\end{gather*}
Here $\vec{x}_{q_i} = x_{q_i}^1\ldots x_{q_i}^{l_i}$ for $x_{q_i}^j\in X$ 
and $\tilde{x}$ denotes $x_{\alpha,t-1}$ for $x = x_{\alpha,t}\in X$, i.e., $R(\tilde{x}) = x$.

In~\eqref{almostRB},~\eqref{LongRel1}, we have 
\begin{gather*}
s\geq1,\quad k,r\geq 0,\quad 
a,b\in RS(X),\quad
z_2,\ldots,z_s\in RS(X)\setminus S(X), \\
\vec{x}_{q_1},\ldots,\vec{x}_{q_{s-1}}\in S(X\setminus X_0),\quad 
\vec{x}_{q_s}\in S(X\setminus X_0)\cup \emptyset
\end{gather*}
and $x_{\beta,r}$ is greater than any letter from $\vec{x}_{q_s}$.

By $R(z_1)$ we denote either that 
$z_1\in RS(X)\setminus S(X)$ or that $R(z_1)$ is absent, i.e., $R(z_1) = \emptyset$.
The same holds for $R(z_{s+1})$.
In particular, the values 
$s=1$, $k=0$, $R(z_1) = R(z_2) = \vec{x}_{q_1} = \emptyset$, 
transform \eqref{LongRel1} to the relation $R(x_{\beta,r}) - x_{\beta,r+1}$.

{\bf Remark 1}. 
In~\eqref{LongRel1}, we use associative words 
$\vec{x}_\alpha\in S(X)$ instead of ordered polynomials from $\Bbbk[X]$, 
otherwise we will have to reduce the products of such polynomials from 
$\Bbbk[X]$ to the ordered ones in all possible compositions from~$S$. 

{\bf Theorem 4}.
The set $S$ is a a Gr\"{o}bner---Shirshov basis in 
$R\As\langle X\rangle$.

{\sc Proof}. 
All compositions between two elements from \eqref{UnivRel} are trivial, 
as it is the method to construct the universal enveloping associative algebra for a given Lie algebra.
Also, compositions of intersection between~\eqref{almostRB} and~\eqref{almostRB} are trivial,
it is a way to get the free associative RB-algebra. 
Thus, all compositions of intersection which are not 
at the same time compositions of inclusion are trivial.

Let us compute a composition of inclusion between \eqref{UnivRel} and \eqref{LongRel1}. Let
\begin{equation}\label{LongWord}
w = R(R(z_1)\vec{x}_{q_1}R(z_2)\ldots R(z_s)\vec{x}_{q_s}x_{\beta,r}x_{\beta,r+1}^k R(z_{s+1}))
\end{equation}
satisfy all conditions described above. 
We apply the relation~\eqref{UnivRel}: $xy-yx-[x,y]$ to the subword
$$
\vec{x}_{q_j} = \vec{x}_{q_j}'xy\vec{x}_{q_j}'',\quad x>y,\ 1\leq j\leq s.
$$ 

Suppose that $j<s$. 
Since the image of an RB-operator is a subalgebra, 
$[x,y]$ lies in the linear combination of elements from $S(X\setminus X_0)$.
Define $w' = w|_{xy\to yx}$ and $w'' = w|_{xy\to [x,y]}$,
here $w|_{\alpha\to\beta}$ means the word $w$ with the subword $\alpha$
replaced by $\beta$.

On the one hand, we have modulo $(S,w)$
\begin{multline*} 
w 
 \mathop{\equiv}\limits^{\eqref{UnivRel}}
 w' + w'' \\
 \mathop{\equiv}\limits^{\eqref{LongRel1}} 
\frac{1}{k+1}\big(
 R(z_1)\ldots R(z_j)\vec{x}_{q_j}|_{xy\to yx+[x,y]}R(z_{j+1}) \ldots R(z_s)\vec{x}_{q_s}x_{\beta,r+1}^{k+1} R(z_{s+1}) \\
 + R(\Sigma_{w'}) + R(\Sigma_{w''})\big),  
\end{multline*}
where $\Sigma_{w}$ is the expression in the RHS in the brackets under the action of~$R$ from~\eqref{LongRB-conseq2}.

On the other hand, modulo $(S,w)$
\begin{multline*}
w \mathop{\equiv}\limits^{\eqref{LongRel1}} 
 \frac{1}{k+1}\big(
 R(z_1)\ldots R(z_j)\vec{x}_{q_j}R(z_{j+1})\ldots R(z_s)\vec{x}_{q_s}x_{\beta,r+1}^{k+1}R(z_{s+1}) + R(\Sigma_{w})\big) \\
 {\mathop{\equiv}\limits^{\eqref{UnivRel}}}
 \frac{1}{k+1}\big(
 R(z_1)\ldots R(z_j)\vec{x}_{q_j}|_{xy\to yx+[x,y]}R(z_{j+1}) \ldots R(z_s)\vec{x}_{q_s}x_{\beta,r+1}^{k+1} R(z_{s+1})
 + R(\Sigma_{w}) \big).
\end{multline*}
So, the composition equals $\dfrac{1}{k+1}(R(\Sigma_{w}) - R(\Sigma_{w'}) - R(\Sigma_{w''}))$
and we may rewrite it briefly as 
$$
\frac{1}{k+1}\sum_{s\in S} R(\alpha_s ([\tilde{x},y]+[x,\tilde{y}]+\lambda[\tilde{x},\tilde{y}] - R^{-1}([x,y]))\beta_s)
$$
for corresponding index set $S$ and $\alpha_s,\beta_s\in RS(X)$.
We get zero, since 
$$
[x,y] = [R(\tilde{x}),R(\tilde{y})] = R([\tilde{x},y]+[x,\tilde{y}]+\lambda[\tilde{x},\tilde{y}])
$$
and $R$ has zero kernel on $L$.

Consider the case $j = s$. The triviality of the corresponding composition of inclusion 
one can derive from the following fact. 
Denote as $A(z_1,\vec{x}_{q_1},z_2,\ldots,z_s,\vec{x}_{q_s},z_{s+1})$ 
the expression~\eqref{LongRB} for $t = l_1+\ldots +l_s+s+1$,
\begin{gather*}
b_1 = z_1,\
b_2 = {x_{q_1}^1}',\
b_3 = {x_{q_1}^2}',\ \ldots,\
b_{l_1+1} = {x_{q_1}^{l_1}}',\
b_{l_1+2} = z_2,\
b_{l_1+3} = {x_{q_2}^1}',\
\ldots, \\
b_{l_1+\ldots +l_{s-1}+s} = z_s,\
b_{l_1+\ldots +l_{s-1}+s+1} = {x_{q_l}^1}',\
\ldots,\ 
b_{l_1+\ldots +l_s+s} = {x_{q_s}^{l_s}}',\
b_{l_1+\ldots +l_s+s+1} = z_{s+1},
\end{gather*}
where $z_1,\ldots,z_{s+1}$ and the words $\vec{x}_{q_i} = x_{q_i}^1 \ldots x_{q_i}^{l_i}$
satisfy the conditions for~\eqref{LongRel1} except the one concerned $x_{\beta,r}$.
We also have that the word 
$\vec{x}_{q_s} = x_{q_1}^1 \ldots x_{q_1}^{l_s}$ contains the biggest letter 
$x_{\beta,r}$, $r\geq1$, on the positions 
$$
K = \{k_1,k_2,\ldots,k_p\mid k_1<k_2<\ldots<k_p\}\subset\{1,2,\ldots,l_s\}.
$$

So, the composition of inclusion is trivial if 
\begin{equation}\label{InductionForComposit}
A(z_1,\vec{x}_{q_1},z_2,\ldots,z_s,\vec{x}_{q_s},z_{s+1}) \equiv 0 \mod (S,w)
\end{equation} 
for $w$ greater than all terns involved in $A(\ldots)$.

To prove~\eqref{InductionForComposit}, we will proceed on by induction on $l_s = |\vec{x}_{q_s}|$. 
For $l_s = 1$, we are done by~\eqref{LongRel1}.

Consider the equality
\begin{equation}\label{GeneralXConvert}
\vec{x}_{q_s}
 = \vec{x}_{q_s,0}x_{\beta,r}^p + \sum\limits_{i=1}^p x_{q_s}^1 x_{q_s}^2\ldots x_{q_s}^{k_i-1}\ldots [x_{q_s}^{k_i},w_i]x_{\beta,r}^{p-i},
\end{equation}
where $\vec{x}_{q_s,0}$ is obtained from $\vec{x}_{q_s}$ by arising all letters $x_{\beta,r}$
with preserving order of all remaining letters,
$w_i$ is obtained from the word $x_{q_s}^{k_i+1}\ldots x_{q_s}^{l_s}$ by arising all $p-i$ letters $x_{\beta,r}$.
For $w_i = w_i^1 w_i^2\ldots w_i^{l_s-k_i}$, $w_i^j\in X$, the bracket $[x_{q_s}^{k_i},w_i]$ in~\eqref{GeneralXConvert} means
$$
[x_{q_s}^{k_i},w_i] 
 = \sum\limits_{j=1}^{l_s-k_i} w_i^1\ldots w_i^{j-1}[x_{q_s}^{k_i},w_i^j]w_i^{j+1}\ldots w_i^{l_s-k_i}.
$$

By~\eqref{LongRel1} and Lemma, we deduce that 
$$
A(z_1,\vec{x}_{q_1},z_2,\ldots,z_s,\vec{x}_{q_s},z_{s+1}) 
 \equiv \sum\limits_{i=1}^p \sum\limits_{j=1}^{l_s-k_i}A(z_1,\vec{x}_{q_1},z_2,\ldots,z_s,\vec{x}_{q_s}(i,j),z_{s+1}) \mod (S,w), 
$$
where 
$\vec{x}_{q_s}(i,j) = x_{q_s}^1\ldots x_{q_s}^{k_i-1}\ldots w_i^1\ldots w_i^{j-1}[x_{q_s}^{k_i},w_i^j]w_i^{j+1}\ldots w_i^{l_s-k_i}x_{\beta,r}^{p-i}$.
The equality 
\linebreak
$A(z_1,\ldots,z_s,\vec{x}_{q_s}(q,i),z_{s+1})\equiv 0$ modulo $(S,w)$ 
follows from the inductive hypothesis.

Consider a composition of inclusion between \eqref{almostRB} and \eqref{LongRel1}.
Let $w = R(w_0)$ be defined by~\eqref{LongWord} and $b\in RS(X)$. 
At first, we have modulo $(S,w)$
\begin{multline*}
R(w_0)R(b) 
 \mathop{\equiv}\limits^{\eqref{almostRB}} 
 R(R(w_0)b+w_0R(b)+\lambda w_0b) \\
 \mathop{\equiv}\limits^{\eqref{LongRel1}} 
 \frac{1}{k+1}  R\big(
 R(z_1)\vec{x}_{q_1}R(z_2)\ldots R(z_s)\vec{x}_{q_s}x_{\beta,r+1}^{k+1} R(z_{s+1})b + R(\Sigma_w)b \big) \\
 + R\big( R(z_1)\ldots R(z_s)\vec{x}_{q_s}x_{\beta,r}x_{\beta,r+1}^k R(z_{s+1})R(b)
 + \lambda R(z_1)\ldots R(z_s)\vec{x}_{q_s}x_{\beta,r}x_{\beta,r+1}^k R(z_{s+1})b \big) \\
 \mathop{\equiv}\limits^{\eqref{LongRel1}}
 \lambda R\big( R(z_1)\ldots R(z_s)\vec{x}_{q_s}x_{\beta,r}x_{\beta,r+1}^k R(z_{s+1})b \big) \\
 + \frac{1}{k+1} \big( R\big( R(z_1)\ldots R(z_s)\vec{x}_{q_s}x_{\beta,r+1}^{k+1} R(z_{s+1})b 
  + R(\Sigma_w)b \big) \\
  + R(z_1)\ldots R(z_s)\vec{x}_{q_s}x_{\beta,r}x_{\beta,r+1}^k R(z_{s+1})R(b) + R(\Sigma_{w|_{R(z_{s+1})\to R(z_{s+1})R(b)}}) \big),
\end{multline*} 
where $\Sigma_{w}$ is the expression in the RHS in the brackets under the action of~$R$ from~\eqref{LongRB-conseq2}.	
At second, modulo $(S,w)$
\begin{multline*}
R(w_0)R(b) 
 \mathop{\equiv}\limits^{\eqref{LongRel1}} 
 \frac{1}{k+1} \big( R(z_1)\vec{x}_{q_1}R(z_2)\ldots R(z_s)\vec{x}_{q_s}x_{\beta,r+1}^{k+1} R(z_{s+1})R(b) + R(\Sigma_w)R(b) \big)\\
 \mathop{\equiv}\limits^{\eqref{almostRB}} 
 \frac{1}{k+1} \big( R(z_1)\ldots R(z_s)\vec{x}_{q_s}x_{\beta,r+1}^{k+1} R(z_{s+1})R(b) 
 + R(R(\Sigma_w)b+\Sigma_wR(b)+\lambda\Sigma_w b) \big).
\end{multline*}
So, the composition of inclusion multiplied by $(k+1)$ equals 
\begin{multline*}
 u = R\big(
 \lambda(k+1) R(z_1)\ldots R(z_s)\vec{x}_{q_s}x_{\beta,r+1}^{k+1}R(z_{s+1})b
 + R(z_1)\ldots R(z_s)\vec{x}_{q_s}x_{\beta,r+1}^{k+1}R(z_{s+1})b \big) \\
 + R(\Sigma_{w|_{R(z_{s+1})\to R(z_{s+1})R(b)}}) - R(\Sigma_wR(b)+\lambda\Sigma_w b).
\end{multline*}
Depending on the last factor, $\Sigma_w$ splits into the sum $\Sigma_w'R(z_{s+1}) + \Sigma_w''z_{s+1}$.
Applying the formulas
\begin{gather*}
\lambda\Sigma_w b = \lambda\Sigma_w' R(z_{s+1})b + \lambda\Sigma_w'' z_{s+1}b, \\
\Sigma_wR(b) = \Sigma_w'R(z_{s+1})R(b) + \Sigma_w''z_{s+1}R(b), \\
\Sigma_{w|_{R(z_{s+1})\to R(z_{s+1})R(b)}}
 = \Sigma_w'R(z_{s+1})R(b) + \Sigma_w''(R(z_{s+1})b+z_{s+1}R(b)+\lambda z_{s+1}b),
\end{gather*}
we deduce
\begin{multline*}
 u = R\big(
 R(z_1)\ldots R(z_s)\vec{x}_{q_s}x_{\beta,r+1}^{k+1}R(z_{s+1})b \\
 + \lambda\big( (k+1)R(z_1)\ldots R(z_s)\vec{x}_{q_s}x_{\beta,r}x_{\beta,r+1}^k R(z_{s+1})b - \Sigma_w' R(z_{s+1})b \big) 
 + \Sigma_w''R(z_{s+1})b \big).
\end{multline*}
Writing down the sum $R(\Sigma_w''R(z_{s+1})b)$, $u$ equals to zero
by the definition of $\Sigma'_w$, $\Sigma''_w$ and Lemma.

Compute a composition of inclusion between \eqref{LongRel1} and \eqref{LongRel1}.
Suppose that we have $w$ defined by~\eqref{LongWord} and 
$$
z_m = R(a_1)\vec{x}_{t_1}R(a_2)\ldots R(a_n)\vec{x}_{t_n}x_{\gamma,o}x_{\gamma,o+1}^m R(a_{n+1})
$$
satisfying all conditions written above for~\eqref{UnivRel}. 

Consider the case $m\leq s$. 
By $\Sigma_{w}$, as earlier, we mean the expression in the RHS in the brackets under the action of~$R$ from~\eqref{LongRB-conseq2} for~$w$.
By $\Delta$, denote the expression in the RHS in the brackets under the action of~$R$ from~\eqref{LongRB-conseq2} for~$z_m$.
We also define 
\begin{gather*}\allowdisplaybreaks
z_m' = R(a_1)\vec{x}_{t_1}R(a_2)\ldots R(a_n)\vec{x}_{t_n}x_{\gamma,o+1}^{m+1} R(a_{n+1}), \\
\widetilde{\Delta} 
 = \Delta - R(a_1)\ldots R(a_n)\vec{x}_{t_n} H R(a_{n+1}), \\
H {=} \sum\limits_{i=2}^{m+1}(-1)^i\!\binom{m+1}{i}\big[x_{\gamma,o},x_{\gamma,o+1}^{(i-1)}\big]x_{\gamma,o+1}^{m+1-i} 
  {+} x_{\gamma,o}x_{\gamma,o+1}^m + x_{\gamma,o+1}x_{\gamma,o}x_{\gamma,o+1}^{m-1} {+} \ldots {+} x_{\gamma,o+1}^m x_{\gamma,o},
\end{gather*}
$\Sigma_w = \Sigma_w' +  \Sigma_w''$,
where we collect all summands from $\Sigma_w$ with the factor $R(z_m)$ in the sum $\Sigma_w'$
and all others in $\Sigma_w''$.

On the one hand, modulo $(S,w)$ we get
\begin{multline} \label{LongRel-LongRelOut}
R\big(R(z_1)\vec{x}_{q_1}R(z_2)\ldots R(z_s)\vec{x}_{q_s}x_{\beta,r} x_{\beta,r+1}^k R(z_{s+1})\big) \\
 \mathop{\equiv}\limits^{\eqref{LongRel1},\ \mbox{out}}  
 \frac{1}{k+1}\big(R(z_1)\vec{x}_{q_1}R(z_2)\ldots \vec{x}_{q_{m-1}} R(z_m)\vec{x}_{q_m} \ldots R(z_s)\vec{x}_{q_s}x_{\beta,r+1}^{k+1}R(z_{s+1}) + R(\Sigma_w)\big) \\
 \mathop{\equiv}\limits^{\eqref{LongRel1},\ \mbox{in}}    
 \frac{1}{(k+1)(m+1)}\big(R(z_1)\vec{x}_{q_1}R(z_2)\ldots \vec{x}_{q_{m-1}}z_m'\vec{x}_{q_m}\ldots R(z_s)\vec{x}_{q_s}x_{\beta,r+1}^{k+1} R(z_{s+1}) \\ 
 + R(z_1)\vec{x}_{\alpha_1}R(z_2)\ldots \vec{x}_{q_{m-1}}R(\Delta)\vec{x}_{q_m}\ldots R(z_s)\vec{x}_{q_s}x_{\beta,r+1}^{k+1} R(z_{s+1})
 + R(\Sigma'_w |_{R(z_m) \to z_m'}) \\
 + R(\Sigma'_w |_{R(z_m)\to R(\Delta)})
 + (m+1)R(\Sigma_w'').
\end{multline}

On the other hand, 
\begin{multline}\label{LongRel-LongRelIn} 
R\big(R(z_1)\vec{x}_{q_1}R(z_2)\ldots R(z_s)\vec{x}_{q_s}x_{\beta,r} x_{\beta,r+1}^k R(z_{s+1})\big) \\ \allowdisplaybreaks
 \mathop{\equiv}\limits^{\eqref{LongRel1},\ \mbox{in}}  
 \frac{1}{m+1}R\big(R(z_1)\vec{x}_{q_1}R(z_2)\ldots \vec{x}_{q_{m-1}} z_m'\vec{x}_{q_m} \ldots R(z_s)\vec{x}_{q_s}x_{\beta,r} x_{\beta.r+1}^k R(z_{s+1}) \\
 + R\big(R(z_1)\vec{x}_{\alpha_1}R(z_2)\ldots \vec{x}_{q_{m-1}} R(\Delta)\vec{x}_{q_m} \ldots R(z_s)\vec{x}_{q_s}x_{\beta,r} x_{\beta,r+1}^k R(z_{s+1})\big) \big) \\
 \mathop{\equiv}\limits^{\eqref{LongRel1},\ \mbox{out}}    
 \frac{1}{(k+1)(m+1)}\big(R(z_1)\vec{x}_{\alpha_1}R(z_2)\ldots \vec{x}_{\alpha_{m-1}}z_m'\vec{x}_{\alpha_m}\ldots R(z_s)\vec{x}_{\alpha_s}x_{\beta,r+1}^{k+1} R(z_{s+1}) \\ 
  + R(\Sigma'_w|_{R(z_m) \to z_m'}) 
  - R(\Sigma_w''|_{z_m\to \widetilde{\Delta}}) \\
  + R(z_1)\ldots \vec{x}_{\alpha_{m-1}}R(\Delta)\vec{x}_{\alpha_m}\ldots R(z_s)\vec{x}_{\alpha_s}x_{\beta,r+1}^{k+1} R(z_{s+1}) + R(\Sigma_w|_{z_m\to \Delta})\big)
   \!\!\mod (S,w).
\end{multline}

Subtracting~\eqref{LongRel-LongRelIn} from~\eqref{LongRel-LongRelOut}, we get 
$\frac{1}{(k+1)(m+1)}R(\Sigma_w''|_{z_m\to C})$, where
$$
C = \Delta - \widetilde{\Delta} - (m+1)z_m = R(a_1)\ldots R(a_n)\vec{x}_{t_r}DR(a_{n+1})
$$
for $D = H - (m+1)x_{\gamma,o}x_{\gamma,o+1}^m$.
Equality of $D$ to zero follows from Lemma.

The proof in the case~$m = s+1$ is slightly different if only $R(a_1) = \emptyset$ and $n = 1$.
Then we need to apply the equality 
$A\big(z_1,\vec{x}_{q_1},z_2,\ldots,z_s,\vec{x}_{q_s}x_{\beta,r}x_{\beta,r+1}^k \vec{x}_{t_1}x_{\gamma,o+1}^{m+1},a_2\big)\equiv 0$
modulo $(S,w)$. Such application is correct, since all terms involved in it have less $R$-degree than~$w$.

It is easy to verify that the remaining compositions of inclusion between~\eqref{UnivRel} and~\eqref{almostRB} are trivial. \hfill $\Square$

{\bf Corollary 1}.
The quotient $A$ of $R\As\langle X\rangle$ by $Id(S)$
is the universal enveloping associative RB-algebra for the Lie algebra 
$L$ with the RB-operator $R$.
Moreover, $L$ injectively embeds into $A^{(-)}$.

{\sc Proof}. 
By~\eqref{almostRB}, $A$ is an associative RB-algebra.
By \eqref{UnivRel}--\eqref{LongRel1}, 
we have that $A$ is enveloping of $L$
for both: the Lie bracket $[,]$ and the action of $R$.
Thus, $A$ is an associative enveloping of $L$.

Let us prove that $A$ is the universal enveloping one.
At first, $A$ is generated by~$L$. At second,
all elements from $S$ are identities in the universal enveloping
associative RB-algebra $U_{RB}(L)$. Indeed,
\eqref{UnivRel} are enveloping conditions for the product,
\eqref{almostRB} is the RB-identity,
the relations~\eqref{LongRel1} as the application of~\eqref{LongRB-conseq2} 
are direct consequences of the RB-identity.

By Theorems 3 and 4 we get the injectivity of embedding $L$ into $A^{(-)}$.
\hfill $\Square$

Let $C$ be a pre- or post-Lie algebra with a linear basis $Y$.
By Theorem~2, $C$ can be injectively embedded into the Lie algebra $L$ with the RB-operator 
$P$ of weight~$\lambda$ and such subset 
$X_0 = \{x_\alpha\mid \alpha\in \Omega\}\subset L$ that
$X = \{x_{\alpha,k}:= P^k(x_\alpha)\mid k\in\mathbb{N},\,\alpha\in\Omega\}$ is a linear basis of $L$.
Here $\Omega$ is a well-ordered set.
Then, by Corollary~1, we embed the Lie RB-algebra $L$ into its universal enveloping associative 
algebra $A$ with the RB-operator~$R$. 
Thus, the subalgebra (in pre- or postalgebra sense) $U(C)$ in $A^{(R)}_\lambda$ 
generated by the set~$Y$ is an enveloping pre- or postassociative algebra of $C$.

Now, we state the main result of the work, the analogue of the Poincar\'{e}---Birkhoff---Witt theorem
for pre-Lie and postLie algebras.

{\bf Theorem 5}.
a) Let $C$ be a pre- or post-Lie algebra, then $U(C)$ is the universal enveloping pre- or postassociative algebra of $C$.

b)\,\cite{Dots} The pairs $(\pre\Lie,\pre\As)$ and $(\post\Lie,\post\As)$ are PBW-pairs.

{\sc Proof}. 
a) Consider the postLie algebra case, the proof when $C$ is a pre-Lie algebra is the same.
Let $Y$ be a basis of $C$.
It is easy to show that $A$ is the universal enveloping associative RB-algebra for $C$
in the sense of the equalities~\eqref{pre/postLie2RBAs}. 
Indeed, given an enveloping associative RB-algebra $D$ of $C$,
define $M$ as the Lie RB-subalgebra of $D^{(-)}$ generated by the image of $C$.
By the universality of~$L$, we have that $M$ is the homomorphic image of the Lie RB-algebra $L$.
Thus, $D$ as the enveloping of $M$ is the homomorphic image of the universal enveloping associative RB-algebra of $M$.
The last one is the homomorphic image of~$A$.

Consider the universal enveloping postassociative algebra $V(C)$ of $C$.
Due to Theorem~2, we may embed $V(C)$ in its universal 
enveloping associative RB-algebra $Z$ with RB-operator $Q$. 
Since $Z$ is the homomorphic image of $A$, 
$V(C)$ as the subalgebra in $Z^{(Q)}_\lambda$ generated by $Y$ is the homomorphic image of $U(C)$.

b) We get it by the construction. \hfill $\Square$

{\bf Corollary 2}. 
The pairs $(\pre\Lie,\RB\As)$ and $(\post\Lie,\RB\As)$ are PBW-pairs.

{\bf Corollary 3}. 
The universal enveloping associative RB-algebra of $U(C)$ is isomorphic to $A$.

As another corollary, we obtain the commutative diagram from Introduction.

\section{Universal enveloping pre/post-associative algebra}

\subsection{Post-Lie case}

Given a postLie algebra $C$ with a linear basis $Y$, 
we want to construct a linear basis of the universal 
enveloping postassociative algebra $U(C)$.
Define $Y^+ = Y\cup Y'\subset X\subset L$ for $Y' = \{P(y)\mid y\in Y\}$.

Given a well-ordered set $Z$, define 
$\Com(Z) := \{w = w_1\ldots w_k\in S(Z)\mid w_j\in Z,\,w_1\leq w_2\leq \ldots \leq w_k\}$.

Let us define a set $E\subset RS(Y)$ by induction. 
At first, $Y\subset E$. At second, the word 
$$
a = R(z_1)w_1R(z_2)\ldots R(z_s)w_sR(z_{s+1})
$$
lies in $E$ for $s\geq1$, $w_1,\ldots,w_s\in \Com(Y^+)$,
$z_2,\ldots,z_s\in E\cap RS(Y)\setminus Y$,
$z_1\in E\cap RS(Y)\setminus Y$
or $R(z_1) = \emptyset$ (the same holds for $z_{s+1}$),
if 

1) at least one of $w_i$ contains a letter from $Y$, 

2) every $z_i$ is not of the following form
\begin{equation}\label{ExceptionForm}
R(q_1)u_1R(q_2)\ldots R(q_t)u_t yR(y)^k R(q_{t+1})
\end{equation}
with the same conditions written for $a$. 
Moreover, $u_1,\ldots,u_t\in \Com(Y')$ and 
$y\in Y$ is greater than all letters from $u_t$. 

{\bf Theorem 6}.
The set $\{e + Id(S) \mid e\in E\}$
forms a linear basis of $U(C)$.

{\sc Proof}.
At first, Theorems 3 and 4 imply that the set $\{e + Id(S) \mid e\in E\}$ is a linearly 
independent set of elements in $R\As\langle Y\rangle$.

At second, we show that $U(C)\subset \Bbbk E + Id(S)$.
Let us prove it by induction of the summary $R$-degree $r$ of factors $a,b$ 
involved in the process of generating the algebra $U(C)$.
By the definition, $Y\subset E$. If $r = 0$, we get only 
linear combinations of elements of the form $w\in S(Y^+)$
with at least one letter from $Y$,
where we identify $R(y)$ with $P(y)$ for every $y\in Y$.
Because of~\eqref{UnivRel}, we have $w\in E + Id(S)$.

Let us prove the inductive step for $r>0$.
Suppose that $a,b\in E$, more precisely, 
$$
a = R(z_1)w_1R(z_2)\ldots R(z_s)w_sR(z_{s+1}),\quad
b = R(q_1)u_1R(q_2)\ldots R(q_t)u_tR(q_{t+1}),\quad
$$
where 
$$
s,t\geq1,\quad 
w_1,\ldots,w_s,u_1,\ldots,u_t\in S(Y^+),\quad
z_2,\ldots,z_s,q_2,\ldots,q_t\in E\cap RS(Y)\setminus Y,
$$
by $R(z_1)$, as above, we mean that 
either $z_1\in E\cap RS(Y)\setminus Y$ or $R(z_1) = \emptyset$.
The same holds for $R(z_{s+1}), R(q_1), R(q_{t+1})$.

To prove Theorem, it is enough to state that $a\succ b,a\prec b,a\cdot b\in \Bbbk E + Id(S)$.
We have modulo $Id(S)$
$$
a\cdot b 
 = \begin{cases}
R(z_1)w_1\ldots R(z_s)w_su_1 R(q_2)\ldots u_tR(q_{t+1}), & R(z_{s+1}) = R(q_1) = \emptyset, \\
R(z_1)w_1\ldots R(z_s)w_s R(z_{s+1}\circ q_1)u_1R(q_2)\ldots u_tR(q_{t+1}), & R(z_{s+1}),R(q_1)\neq \emptyset, \\
ab, & \mbox{otherwise},
\end{cases}
$$
where $z_{s+1}\circ q_1 = z_{s+1}\succ q_1 + z_{s+1}\prec q_1 + z_{s+1}\cdot q_1$.
In the third case, we have an element from $E$. In the first one, we need to express 
$w_su_1\in U(L)$ via basic elements from $\Com(Y^+)$ and so $a\cdot b\in \Bbbk E$.
Finally, in the second case, after rewriting $z_{s+1}\circ q_1$ as a linear combination 
of elements from $E + Id(S)$, we only need to check the condition~2 of the definition of $E$.
If it's required, we apply the relation~\eqref{LongRel1} and we are done by induction on $r$.

We have modulo $Id(S)$
$$
a\succ b 
 = \begin{cases}
R(a)u_1 R(q_2)\ldots u_tR(q_{t+1}), & R(q_1) = \emptyset, \\
R(a\circ q_1)u_1R(q_2)\ldots u_tR(q_{t+1}), & R(q_1)\neq \emptyset. \\
\end{cases}
$$
By the inductive hypothesis, $a\circ q_1\in \Bbbk E + Id(S)$.
We need to check the condition~2 of the definition of $E$
for the first $R$-letter of the product. 
As above, we apply~\eqref{LongRel1}.

The case of $a\prec b$ can be considered analogously. \hfill $\Square$

{\bf Remark 2}.
We may reformulate Theorem 6 in terms avoiding $Id(S)$
and any factorization. For this, we need to define by induction
the products $\succ,\prec,\cdot$ on the space~$\Bbbk E$. 

\subsection{Pre-Lie case}

Given a pre-Lie algebra $C$ with a linear basis $Y$, 
we construct a linear basis of the universal enveloping 
preassociative algebra $U(C)$.

Let us define a set $E\subset RS(Y)$ by induction. 
At first, $Y\subset E$. At second, the word 
$$
a = R(z_1)w_1R(z_2)\ldots R(z_s)w_sR(z_{s+1})
$$
lies in $E$ for $s\geq1$, $w_1,\ldots,w_s\in \Com(Y^+)$,
$z_2,\ldots,z_s\in E\cap RS(Y)\setminus Y$,
$z_1\in E\cap RS(Y)\setminus Y$
or $R(z_1) = \emptyset$ (the same holds for $z_{s+1}$),
if 

1) the only $w_i$ contains a letter from $Y$, 

2) every $z_i$ is not of the following form
\begin{equation}\label{ExceptionForm2}
R(q_1)u_1R(q_2)\ldots R(q_t)u_t yR(y)^k R(q_{t+1})
\end{equation}
with the same conditions written for $a$; 
$u_1,\ldots,u_t\in \Com(Y')$ and 
$y\in Y$ is greater than all letters from $u_t$. 

{\bf Theorem 7}.
The set $\{e + Id(S) \mid e\in E\}$
forms a linear basis of $U(C)$.

{\sc Proof}.
Analogously to the proof of Theorem~6. \hfill $\Square$

{\bf Example}.
If $C$ is a one-dimensional pre-Lie algebra with linear basis $\{y\}$,
then $\{y(P(y))^k\}$ is a linear basis of $U(C)$. 
Indeed, because of the condition~2, all basic elements have no $R$-letters,
so we have unique $y$ and several $P(y)$ to form the elements from~$E$.
Moreover, in the case $C = \Span\{y\}$ we have the isomorphism 
$U_{gr}(C)\cong \pre\Com\langle y\rangle$,
where $U_{gr}(C)$ means the associated graded preassociative algebra 
obtained from $U(C)$ by the filtration by the degree
and $\pre\Com\langle Z\rangle$ means a free precommutative algebra generated by the set $Z$.
If $C = \Span\{y\}$ has the trivial product, then $U(C)\cong\pre\Com\langle y\rangle$.

{\bf Remark 3}. 
Note that given a pre-Lie algebra $C$, 
the injective enveloping preassociative algebra $T$ constructed in~\cite{Gub2018-pre}
is isomorphic to $U(C)$, so it is the universal enveloping one.
In the current paper, we have embedded $C$ into its universal enveloping Lie RB-algebra,
but in~\cite{Gub2018-pre} $C$ was embedded into the enveloping Lie RB-algebra $\hat{C}$
arisen from the proof of Theorem~2. A posteriori, we conclude that 
it is enough to embed $C$ into the doubling (as a vector space) Lie RB-algebra
to preserve all required connections for the construction of the universal 
enveloping preassociative algebra $U(C)$.

\section*{Acknowledgements}

This work was supported by the Austrian Science Foundation FWF grant P28079.

\noindent Vsevolod Gubarev \\
University of Vienna \\
Oskar-Morgenstern-Platz 1, 1090, Vienna, Austria \\
e-mail: vsevolod.gubarev@univie.ac.at

\begin{thebibliography}{99}
\bibitem{Aguiar00}
M. Aguiar, Pre-Poisson algebras, {\it Lett. Math. Phys.} {\bf 54} (2000) 263--277.

\bibitem{BBGN2012}
C. Bai, O. Bellier, L. Guo, X. Ni, Splitting of operations, Manin products, and Rota---Baxter operators,
{\it Int. Math. Res. Notices} {\bf 3} (2013) 485--524.

\bibitem{BaiGuoNi10}
C. Bai, L. Guo, X. Ni, $\mathcal{O}$-operators on associative algebras, associative
Yang---Baxter equations  and dendriform algebras, Conference Proceedings 
{\it Quantized Algebra and Physics} (2012) 10--51.

\bibitem{Baxter60}
G. Baxter, An analytic problem whose solution follows from a simple algebraic identity, 
{\it Pacific J. Math.} {\bf 10} (1960) 731--742.

\bibitem{BelaDrin82}
A.A. Belavin, V.G. Drinfel'd, Solutions of the classical Yang-Baxter equation for
simple Lie algebras, {\it Funct. Anal. Appl.} {\it 16} (3) (1982) 159--180.

\bibitem{BokutChen}
L.A. Bokut, Yu. Chen, X. Deng, Gr\"{o}bner-Shirshov bases for Rota-Baxter algebras,
{\it Sib. Math. J.} {\bf 51} (6) (2010) 978--988.

\bibitem{Burde06}
D. Burde,
Left-symmetric, or pre-Lie algebras in geometry and physics,
{\it Cent. Eur. J. Math.}, {\bf 4} (3), 325--357 (2006).

\bibitem{Burde12} 
D. Burde, K. Dekimpe and K. Vercammen, Affine actions on Lie groups and post-Lie algebra structures,
{\it Linear Algebra Appl.} {\bf 437} (2012), no. 5, 1250--1263.

\bibitem{Cartier72}
P. Cartier, On the structure of free Baxter algebras, {\it Adv. Math.} {\bf 9} (1972) 253--265.

\bibitem{Chapaton}
F. Chapoton, Un th\'{e}or\`{e}me de Cartier-Milnor-Moore-Quillen pour les big\`{e}bre dendriformes
et les alg\`{e}bre braces, {\it J. Pure Appl. Algebra} {\bf 168} (2002), 1--18.

\bibitem{Chen11}
Yu. Chen, Q. Mo, Embedding dendriform algebra into its universal enveloping Rota---Baxter algebra,
{\it Proc. Amer. Math. Soc.} {\bf 139} (12) (2011) 4207--4216.

\bibitem{Dokas13}
I. Dokas, Pre-Lie algebras in positive characteristic, {\it J. Lie Theory} {\bf 23} (4) (2013) 937--952.

\bibitem{Dots}
V. Dotsenko, P. Tamaroff, Endofunctors and Poincar\'{e}-Birkhoff-Witt theorems, 
arXiv:1804.06485 [math.CT].

\bibitem{Fard02}
K. Ebrahimi-Fard, Loday-type algebras and the Rota-Baxter relation, 
{\it Lett. Math. Phys.} {\bf 61} (2002) 139--147.

\bibitem{FardGuo07}
K. Ebrahimi-Fard, L. Guo, Rota---Baxter algebras and dendriform algebras,
{\it J. Pure Appl. Algebra} {\bf 212} (2) (2008) 320--339.

\bibitem{ELMM} 
K. Ebrahimi-Fard, A. Lundervold, I. Mencattini, H. Z. Munthe-Kaas, 
Post-Lie Algebras and Isospectral Flows, {\it Symmetry Integr. Geom.} {\bf 11} (2015), Paper~093, 16~pp.

\bibitem{Gerst63}
M. Gerstenhaber, The cohomology structure of an associative ring,
{\it Ann. of Math.} {\bf 78} (1963) 267--288.

\bibitem{GinzKapr}
V. Ginzburg, M. Kapranov, Koszul duality for operads, {\it Duke Math. J.}
{\bf 76} (1) (1994) 203--272.

\bibitem{GolubchikSokolov}
I.Z. Golubchik, V.V. Sokolov, Generalized operator Yang-Baxter equations, integrable ODEs
and nonassociative algebras, {\it J. Nonlinear Math. Phys.} {\bf 7} (2) (2000) 184--197.

\bibitem{Gub2017Com}
V. Gubarev,
Universal enveloping commutative Rota---Baxter algebras of pre- and post-commutative algebras,
{\it Axioms} {\bf 6} (4) (2017) 1--33.

\bibitem{Gub2017As}
V. Gubarev,
Universal enveloping associative Rota---Baxter algebras of preassociative
and post\-associative algebras, {\it J. Algebra} {\bf 516} (2018) 298--328.

\bibitem{Gub2017Lie}
V. Gubarev,
Universal enveloping Lie Rota---Baxter algebras of pre-Lie and post-Lie algebras,
{\it Algebra and Logic} (accepted), arXiv:1708.06747.

\bibitem{Gub2018-post}
V. Gubarev, 
Embedding of postLie algebras into postassociative algebras, \\
arXiv:1808.08839 [math.RA].

\bibitem{Gub2018-pre}
V. Gubarev, 
Embedding of pre-Lie algebras into preassociative algebras, 
{\it Algebr. Colloq.} (accepted), arXiv:1808.09822 [math.RA].

\bibitem{GubKol2013}
V. Gubarev, P. Kolesnikov, 
Embedding of dendriform algebras into Rota---Baxter algebras,
{\it Cent. Eur. J. Math.} {\bf 11} (2) (2013) 226--245.

\bibitem{GubKol2014}
V. Gubarev, P. Kolesnikov, Operads of decorated trees and their duals, 
{\it Comment. Math. Univ. Carolin.} {\bf 55} (4) (2014) 421--445.

\bibitem{Guo2011}
L. Guo, {\it An Introduction to Rota---Baxter Algebra}, Intern. Press, Somerville, MA; Higher education press, Beijing, 2012. 

\bibitem{Guo2013}
L. Guo, W. Sit, and R. Zhang, Differential type operators and Gr\"{o}bner-Shirshov bases, 
{\it J. Symb. Comput.} {\bf 52} (2013) 97--123.

\bibitem{Kol2017}
P. Kolesnikov, Gr\"{o}bner---Shirshov bases for pre-associative algebras,
{\it Commun. Algebra} {\bf 45} (12) (2017) 5283--5296.

\bibitem{Koszul61}
J.-L. Koszul,
Domaines born\'{e}s homog\`{e}nes et orbites de groupes de
transformations affines,
{\it Bull. Soc. Math. France}, {\bf 89}, 515--533 (1961).

\bibitem{Loday95}
J.-L. Loday,
Cup-product for Leibniz cohomology and dual Leibniz algebras, {\it Math. Scand.} {\bf 77} (2) (1995) 189--196.

\bibitem{Dialg99}
J.-L. Loday, Dialgebras, {\it Dialgebras and related operads}, Springer-Verl., Berlin, 2001, 1--61. 

\bibitem{Loday2007}
J.-L. Loday, On the algebra of quasi-shuffles, {\it Manuscripta Math.} {\bf 123} (2007) 79--93.

\bibitem{LP93}
J.-L. Loday, T. Pirashvili, Universal enveloping algebras of Leibniz algebras and (co)homology,
{\it Math. Ann.} {\bf 296} (1993) 139--158.

\bibitem{Trialg01}
J.-L. Loday, M. Ronco, Trialgebras and families of polytopes, {\it Comtep. Math.} {\bf 346} (2004) 369--398.

\bibitem{Manchon}
D. Manchon,
A short survey on pre-Lie algebras,
{\it Noncommutative Geometry and Physics: Renormalisation,
Motives, Index Theory}, 89--102 (2011).

\bibitem{PBW}
A.A. Mikhalev, I.P. Shestakov, PBW-pairs of varieties of linear algebras,
{\it Commun. Algebra} {\bf 42} (2) (2014) 667--687.

\bibitem{sl2}
Yu Pan, Q. Liu, C. Bai, L. Guo,
PostLie algebra structures on the Lie algebra $\mathrm{sl}(2,\mathbb{C})$,
{\it Electron. J. Linear Al.}, {\bf 23}, 180--197 (2012).

\bibitem{Ronco}
M. Ronco, Eulerian idempotents and Milnor---Moore theorem for certain non-cocommutative Hopf algebras, 
{\it J. Algebra} {\bf 254} (2002), 152--172.

\bibitem{Rota68}
G.-C. Rota, Baxter algebras and combinatorial identities I, {\it Bull. Amer. Math. Soc.} {\bf 75} (1969) 325--329.

\bibitem{Semenov83}
M.A. Semenov-Tyan-Shanskii, What is a classical $r$-matrix? {\it Funct. Anal. Appl.} {\bf 17} (4) (1983) 259--272.

\bibitem{Vallette2007}
B. Vallette, Homology of generalized partition posets, {\it J. Pure Appl. Algebra} {\bf 208} (2) (2007) 699--725.

\bibitem{Vinberg63}
E.B. Vinberg,
The theory of homogeneous convex cones, {\it Tr. Mosk. Mat. Obs.}, {\bf 12}, 303--358 (1963).
\end{thebibliography}
\end{document}